\documentclass[final,3p,times]{elsarticle}

\usepackage{latexsym}
\usepackage{amssymb}
\usepackage[mathscr]{eucal}
\usepackage{epsfig}
\usepackage{lscape}
\usepackage{amssymb}
\usepackage{amsmath}
\usepackage{graphicx}
\newcommand{\rr}{\mathbb{R}}
\newcommand{\zz}{\mathbb{Z}}
\newcommand{\cc}{\mathbb{C}}
\newcommand{\nn}{\mathbb{N}}
\newcommand{\dis}{\displaystyle}
\newtheorem{Theorem}{Theorem}

\newtheorem{Lemma}[Theorem]{Lemma}
\newtheorem{Corollary}[Theorem]{Corollary}

\newtheorem{Remark}[Theorem]{Remark}

\renewcommand\Re{\operatorname{Re}}

\def\a{\alpha}
\def\b{\beta}

\def\be{\begin{equation}}
\def\ee{\end{equation}}
\def\bea{\begin{eqnarray}}
\def\eea{\end{eqnarray}}
\def\bean{\begin{eqnarray*}}
\def\eean{\end{eqnarray*}}

 \def\proof#1. {\par
                      \ifdim\lastskip<15pt
                      \removelastskip\penalty-200
                      \vskip5pt plus3pt minus3pt
                      \fi
                       {\def\a{#1}
                       \ifx\a\empty
                       {\noindent\bf Proof.}
                       \else
                       {\noindent\bf Proof of #1.}
                       \fi}\enspace}
\def\endproof{\hfill\hspace{-6pt}\rule[-4pt]{6pt}{6pt}
\vskip8pt plus3pt minus 3pt}

\journal{...}

\begin{document}

\begin{frontmatter}


\author{A Jooste\fnref{label3}}
\ead{alta.jooste@up.ac.za}\address{Department of Mathematics and
Applied Mathematics, University of Pretoria, Pretoria, 0002, South
Africa}

\author{K Jordaan\fnref{label2}}
\ead{kjordaan@up.ac.za} \fntext[label2]{Research by these
authors is partially supported by the National Research Foundation
under grant number 2054423.} \address{Department of Mathematics and Applied Mathematics, University
of Pretoria, Pretoria, 0002, South Africa}

\author{F To\'{o}kos}\ead{ferenc.tookos@helmholtz-muenchen.de}\address{Institute for Biomathematics and Biometry, Helmholtz Zentrum M\"{u}nchen, Neuherberg, Germany}

\title{Zeros of Meixner and Krawtchouk polynomials}

\begin{abstract} We investigate the zeros of a family of
hypergeometric polynomials $\phantom{ }_2F_1(-n,-x;a;t)$,
$n\in\nn$ that are known as the Meixner polynomials for certain values of the
parameters $a$ and $t$. When $a=-N$, $N\in\nn$ and $t=\frac1{p}$, the
polynomials $K_n(x;p,N)=(-N)_n\phantom{ }_2F_1(-n,-x;-N;\frac1{p})$, $n=0,1,\dots N$,
$0<p<1$ are referred to as Krawtchouk polynomials. We prove results for the zero location of the orthogonal polynomials $K_{n}(x;p,a)$, $0<p<1$ and $a>n-1$, the quasi-orthogonal polynomials $K_{n}(x;p,a)$, $k-1<a<k$, $k=1,\dots,n-1$ and $p>1$ or $p<0$  as well as the non-orthogonal polynomials $K_{n}(x;p,N)$, $0<p<1$ and $n=N+1,N+2,\dots$.
We also show that the polynomials $K_{n}(x;p,a)$, $a\in \rr$ are real-rooted when $p\rightarrow 0$.
\end{abstract}


\end{frontmatter}

\section{Introduction}

\medskip \noindent
Polynomials that are real-rooted are of interest in various branches
of mathematics, including combinatorics and approximation theory.
Real-rootedness is used to prove the log-concavity and unimodality
of the sequence of coefficients of a polynomial which otherwise
demands cumbersome combinatorial manipulations (see e.g.\
\cite{Stanley89} and \cite{Brenti94}) and is useful in establishing
convergence properties of rational approximants  (cf. \cite{pade}).

\medskip \noindent The real-rootedness of a given polynomial is a classical
problem and various tools are available to prove that the zeros are
real, for example, in \cite{Br} and \cite{PF}, P\'olya frequency
sequences are used to show that some hypergeometric polynomials have
only real roots. One of the standard ways to prove real-rootedness
of polynomials, is to use their relation to orthogonal polynomials.

\medskip \noindent
A sequence of polynomials $\{{P_n}\}_{n=0}^N$, $N\in\nn \cup \{ \infty \} $ is orthogonal
with respect to some positive Borel measure $\phi$, if
\begin{equation}\int_{S}P_n(x)P_m(x)~d\phi(x)= d_n^{~2} ~\delta_{mn},~m,n=0,1,\dots N \label{TTRR}\end{equation}
where $S$ is the support of the measure and $\delta_{mn}$ is Kronecker's symbol. It is a well-known
classical result that the $n$ zeros of $P_n$ are real, distinct and
lie in the convex hull of $\mbox{supp}(\phi)$ (cf. \cite{Szego}).

\medskip \noindent
If the measure $\phi$ is discrete with weight $\rho_i$ at the distinct points
$\{x_i\}_{i=0}^{N-1}$, then (\ref{TTRR}) becomes (cf. \cite{Atk})
\begin{eqnarray*}\sum_{i=0}^{N-1} P_n(x_i)P_m(x_i)~\rho_i~=~ d_n^{~2}~\delta_{mn},~m,n=0,1,\dots N\end{eqnarray*}
and the sequence $\{P_n\}_{n=0}^{N}$ is discrete orthogonal on $(0,N-1)$.

\medskip \noindent The classical orthogonal polynomials of a discrete variable
are used extensively and have many applications in  various fields
including combinatorial analysis, theoretical and mathematical
physics, group representation theory and stochastic processes. For
example, Meixner polynomials have been used to analyse discrete
stochastic processes in the context of spectral analysis in the
Laplace domain (cf. \cite{Graf}) and close relationships have been
found linking generalised spherical harmonics for SU(2) with a
special class of Meixner polynomials known as Krawtchouk polynomials
\cite{Koorn}.

\medskip In this paper we make a comprehensive study of the zero location of Meixner and Krawtchouk polynomials,
in particular for parameter values where (some of) the zeros are real.

\medskip \noindent Meixner polynomials are defined in terms of
the $_2F_1$ hypergeometric function (cf. \cite [p.174]{Ismail})
\bean M_n(x;\beta,c)&=&(\beta)_n\phantom{ }~_2F_1(-n,-x;\beta;1-\frac 1c)\\\
&=&(\beta)_n\phantom{ }~
\sum_{k=0}^{n}\frac{(-n)_k(-x)_k(1-\frac 1c)^k}{(\beta)_kk!},\phantom{ }~ \beta,c\in
\rr, ~\beta\neq-1,-2,\dots,-n+1, ~c\neq 0\eean
where $(~ )_n$ is Pochhammer's symbol defined by
\begin{eqnarray*}(\a)_n& =&(\a)(\a+1).....(\a+n-1)~\mbox{for}~
n\geq1\\
(\a)_0&=&1 ~\mbox{when}~ \a\neq0.\end{eqnarray*}

\medskip \noindent
The sequence $\{M_n(x;\b,c)\}_{n=0}^{\infty}$ satisfies the
discrete orthogonality relation (cf. \cite[p.346]{AAR})
\begin{equation}\label{1}\sum_{x=0}^{\infty}\frac{c^x(\b)_x}{x!}M_m(x;\b,c)M_n(x;\b,c)=
\frac{(\b)_nn!}{c^n(1 - c)^{\b} }\phantom{}~\delta_{mn}\end{equation}
when $0<c<1$ and $\b>0$, and hence the zeros are real, distinct and in $(0,\infty)$ for these values of the
parameters $\b$ and $c$. It can also be shown that, for the same parameter values, the Meixner
polynomials have a
non-standard orthogonality property defined in terms of a discrete inner
 product involving difference operators (cf. \cite{MPPR}).

\medskip \noindent When $c>1$ and $\b>0$ the orthogonality relation (\ref{1}) can be written
as (cf. \cite[p. 177 (3.7)]{Chihara78})
\[\sum_{x=-\infty}^{-\b} \frac{(\b)_{-\b-x}}{(-\b-x)!}\phantom{ }~(c)^{\b+x}M_m(x;\b,c)M_n(x;\b,c)=
\frac{(\b)_n c^\b n!}{c^n (c - 1)^\b}\phantom{ }~\delta_{mn}\] by an
application of Pfaff's transformation (cf. \cite[(2.2.6)]{AAR})
where \[\sum_{x=-\infty}^{-\b}f(x)=f(-\b)+f(-\b-1)+f(-\b-2)+\dots\]
and therefore the zeros of the polynomials
$\{M_n(x;\b,c)\}_{n=1}^{\infty}$ are distinct and in $(-\infty,-\b)$.

\medskip \noindent
Another non-standard orthogonality for Meixner
polynomials $M_n(x;\b,c)$ when $\b,~c \in \cc$, $c \notin [0,\infty)$ and $-\b
\notin \nn$  is discussed in \cite{costas}. The standard orthogonality for a finite
number of Meixner polynomials $M_n(x;\b,c)$ when $c<0$ and $\b$ is
equal to a negative integer, say $\b=-N$, $N\in\nn$ is that of the Krawtchouk polynomials
 defined by (cf.
\cite{Ismail}, p. 183)
\bea\label{Kraw} K_n(x;p,N)&=&(-N)_n\phantom{ }~_2F_1(-n,-x;-N;\frac
1p),~n=0,1,\dots,N.\eea The
orthogonality relation for Krawtchouk polynomials is given by
\[\sum_{x=0}^{N}w(x;p,N) K_m(x;p,N)K_n(x;p,N)=0\] when $m < n \leq
N$; $m,n,N\in \nn$ and $0<p<1$ where $\dis{w(x;p,N)={N \choose
x}\left(p\right)^x\left(1-p\right)^{N-x}}$  is positive at the
mass points $x=0,1,\dots,N$ of the discrete measure. This implies that for $0<p<1$ and $n\leq N$,
 $n,N\in \nn$, the zeros of $K_n(x;p,N)$ are
real, distinct and in the interval $(0,N)$. Furthermore, by an argument due to L. Fej\'er
there is at most one zero of $K_n(x;p,N)$ in between any two
consecutive mass points (cf. \cite{Jordan}) and in the particular case
where $n=N$, the zeros of $K_n(x;p,n)$, denoted by $x_i$,
$i=1,2,\dots,n$ in ascending order, interlace with the mass points
as follows \begin{eqnarray}\label{intmass}0<x_1<1<x_2<2<\dots<x_n<n.\end{eqnarray}

\noindent The zero asymptotics of normalised Krawtchouk polynomials when the ratio of parameter
$n/N \to \alpha$ as $n,\ N \to \infty$ was investigated in \cite{Saff97} and \cite{Saff00} by
finding the support and density of the constrained extremal measure for all possible values of the parameter $\alpha$ and
the asymptotic zero distribution of Meixner polynomials has also been studied by various authors
(cf. \cite{FLS} and \cite{JW}).

\medskip \noindent Since $\displaystyle(-a+k)_{n-k}=\frac{(-a)_n}{(-a)_k}$, (\ref{Kraw}) can be rewritten as
\bea\label{Kraw2} K_n(x;p,a)&=&\sum_{k=0}^n\frac{(-n)_k(-x)_k(-a+k)_{n-k}}{k!p^k}\\ &=&
M_n(x;-a,\frac{p}{p-1})\label{m}\eea  which is valid for every $n$ and
can be used to define Krawtchouk polynomials for any $n\in \nn$ and $a\in\rr$.

\medskip \noindent The polynomials $K_n(x;p,a)$ have the standard orthogonality of Meixner and Krawtchouk polynomials for the parameter ranges \begin{itemize}
\item[] $p<0$, $a<0$,  $n=0,1,2,\dots$;
\item[] $p>1$, $a<0$, $n=0,1,2,\dots$~\mbox{and}
\item[] $0<p<1$, $a=N$, $N\in\nn$, $n=0,1,2\dots N$,
\end{itemize} hence our discussion will focus on the zeros of polynomials $K_n(x;p,a)$ for parameter values
\begin{itemize}
\item[(i)] $0<p<1,$ $a>n-1$,~ $n\in\nn$;
\item[(ii)] $p<0$, $a>0$, $n=0,1,2,\dots$;
\item[(iii)] $p>1$, $a>0$, $n=0,1,2,\dots$;
\item[(iv)] $0<p<1$, $a=N$, ~$n=N+1,N+2,\dots$;
\item[(v)] $0<p<1$, $a<0$ and
\item[(vi)] $p\to 0$, $a\in \rr$, $n=0,1,2\dots$.
\end {itemize}

\medskip \noindent
We begin with case (i) in Section \ref{22} where we extend
the conclusion following from the discrete orthogonality of Krawtchouk polynomials for integer values of the parameter $N$, $N+1>n$, $0<p<1$ to prove that the zeros of Krawtchouk polynomials $K_n(x;p,a)$ are real, distinct and lie in the interval $(0,a)$ for
all real values of the parameter $a$, $a>n-1$ and $0<p<1$. For $a\neq0,1,\dots,n-1$ we obtain \bea
\nonumber K_n(x;p,a)&=&(-a)_n\phantom{ }~_2F_1\left(-n,-x;-a; \frac 1p\right)\\
&=&(-a)_n (1-\frac 1p)^n \phantom{ }~_2F_1\left(-n, -a+x;
-a; \frac {\frac{1}{p}}{\frac{1}{p} -1}\right) ~~~\mbox{(cf. \cite[(2.2.6)]{AAR})}\label{77}\\
&=&(1-\frac 1p)^n \phantom{ }K_n(a-x;1-p,a),\label{symmetry}
\eea
\noindent
a general symmetry property of the Krawtchouk polynomials, since by continuity
it holds for $a\in \rr$, which we will use in the proof of our theorem. Special cases of this symmetry property were proved
for $a=-N$, $N\in \nn$ and $x=0,\dots,N$ (cf. \cite{Ismail98, Chihara90}) using the generating
function for Krawtchouk polynomials (cf. \cite{Koekoek}) and for the Meixner polynomials when
$x=0,\dots, n$ in \cite{Chihara78}.

\medskip\noindent In Section \ref{44} we consider cases (ii) and (iii), proving that the polynomials $K_n(x;p,a)$
are quasi-orthogonal of order $k$ for $k-1<a<k$, $k=1,\dots,n-1$ and $p<0$ or $p>1$.

\medskip \noindent
We discuss case (iv) in Section \ref{3}. Results obtained in \cite{Ask},
\cite{AskSylv}, \cite{Holtz} and \cite{Sylvester} for the zeros of Krawtchouk
polynomials $K_n(x;p,N)$ of degree $n=N+1$, where these polynomials are no longer
orthogonal, are extended  to polynomials of degree $n=N+2$ and $n=N+3$. We make use of the product decomposition
(cf. \cite{AskSylv} and \cite{costas})
\begin{eqnarray}\label{product}K_n(x;p,N)=K_{N+1}(x;p,N)K_{n-N-1}(x-N-1;p,-N-2),~ p\neq0,1,~n>N \in \nn\end{eqnarray} which also shows that for case (iv) it suffices to study polynomials $K_n(x;p,N)$ for $0<p<1$, $N=-2,-3,\dots$ (case (v) for integer values of $a$).

\medskip \noindent In the last section we prove that the polynomials
$\{K_n(x;p,a)\}_{n=1}^{\infty}$ are real-rooted
when $a\in \rr$ and $p\to 0$ (case (vi)).

\medskip \noindent
Observe that for the special case when $p=1$ and $a\in \rr$ we have

\begin{eqnarray*}K_n(x;1,a)&=&(-a)_n\phantom{}_2F_1(-n,-x;-a;1)\\
&=&(x-a)(x-a+1).....(x-a+n-1)\end{eqnarray*}
\noindent which vanishes when $x=a,a-1,\dots,a-n+1$
whereas for $p\to \infty$ the polynomial
\[ (-a)_n\phantom{ }~_2F_1(-n,-x;-a;0)
=(-a)_n\]
 has $n$ zeros at infinity if it is considered as a
polynomial of degree $n$ in $x$ (cf. \cite[(6.72.3)]{Szego}).

\medskip \noindent
The only remaining case is $0<p<1$, $a>0$ where evidence from numerical examples indicates that, in general, there are no real zeros apart from the one real zero when $n$ is odd.

\section{The zeros of $K_{n}(x;p,a)$, $0<p<1$ and $a>n-1$}\label{22}

\noindent The finite sequence
$\{K_n(x;p,N)\}_{n=0}^{N}$ satisfies the three term recurrence
relation (cf. \cite[(1.10.3)]{Koekoek})
\begin{eqnarray}xK_{n}(x;p,N)&=&A_nK_{n+1}(x;p,N)+B_nK_n(x;p,N)+C_nK_{n-1}(x;p,N)\label{ttR}\\ \nonumber K_0&=&1,~~K_{-1}=0\end{eqnarray}
where
$A_n=p$ and $C_n=n(1-p)(N-n+1)$ so that
$A_{n-1}C_n>0$ when $0<p<1$ and $n<N+1$ regardless of
whether $N$ is an integer or not.

\medskip\noindent It follows from a theorem often
attributed to Favard (cf. \cite{Askey}) that there is at least one positive measure
$d\alpha(x)$ so that
\[\int_{-\infty}^{\infty}K_n(x;p,a)K_m(x;p,a)d\alpha(x)=0,~ m\neq n,~a>n-1,~0<p<1\]
and hence $K_n(x;p,a)$ has $n$ real, distinct zeros
when $n<a+1$, $0<p<1$.  However, the set containing the real zeros does
not follow immediately.

\medskip \noindent We prove that the zeros of $K_n(x;p,a)$ for $a>n-1$ and $0<p<1$ are in $(0,a)$ by using
a generalised Sturmian sequence argument applied to
solutions of difference equations (cf. \cite{Porter}), as was done
in \cite{Levit} for Hahn polynomials. We begin by proving that if $r$ denotes a zero of $K_n(x;p,a)$ in $(0,a)$,
then $r-1$ and $r+1$ cannot be zeros of $K_n(x;p,a)$ and, in
addition, there will be an odd number of zeros of $K_n(x;p,a)$ in
the interval $(r-1,r+1)$.

\begin{Lemma}\label{l} Let $a\in\rr$, $n \in\nn$, $a>n-1$ and $0<p<1$. If $r$ is a zero
of $K_n(x;p,a)$ and $r\in(0,a)$, then\newline
$K_n(r-1;p,a)K_n(r+1;p,a)<0$.\end{Lemma}

\proof. Let $a\in\rr$, $a+1>n$ and $0<p<1$. Consider the difference
equation (cf \cite[(1.10.5)]{Koekoek})
\begin{eqnarray}A(x)K_n(x+1;p,a)+C(x)K_n(x-1;p,a)&=&B(x)K_n(x;p,a)\label{2}~\mbox{where}\\
 \nonumber A(x)&=& p(x-a)\\
\nonumber C(x)&=&x(p-1)\\
\nonumber B(x)& =& n+p(x-a)+x(p-1).\end{eqnarray} Note that $A(x)<0$ and
$C(x)<0$ when $x\in(0,a)$ and $0<p<1$.

\medskip \noindent
Suppose $r$ is a zero of $K_n(x;p,a)$ in the interval $(0,a)$,
then
\begin{equation}A(r)K_n(r+1;p,a)+C(r)K_n(r-1;p,a)= 0.\label{51}\end{equation}
Assume that
\begin{equation}K_n(r+1;p,a)=0\label{4}.\end{equation} Letting $x=r+1$ in (\ref{2}) we obtain $A(r+1)K_n(r+2;p,a)=0$
and if $r+1\in(0,a)$,
it follows that
$A(r+1)<0$ and $K_n(r+2;p,a)=0$. By repeating this argument we
can prove that
\begin{equation}K_n(r+i;p,a)=0\text{ for all}\ i \text{ such that }\ 0<r+i-1<a.\label{5}\end{equation} 
Under our assumption (\ref{4}), it also follows from equation
(\ref{51}) that $C(r)K_n(r-1;p,a)=0$ if $r\in(0,a)$
and since $C(r)<0$ for these values of $r$, we have that
$K_n(r-1;p,a)=0 \text{ if }r\in(0,a)$.
In the same way as before we can prove that
\begin{equation}K_n(r-j;p,a)=0\text{ for all}\ j \text{ such that }\ 0<r-j+1<a.\label{6}\end{equation}\ 
In short, it follows from results (\ref{5}) and (\ref{6}) that
$K_n(x;p,a)$ has as zeros all numbers $r+i$, $i\in \zz$ with
$-1<r+i<a+1$. This means that $K_n(x;p,a)$ has a total of at most
$\lfloor a+1-(-1)\rfloor=\lfloor
a+2\rfloor=\lfloor a\rfloor+2\geq n+1>n$ zeros
unless both $a$ and $r$ are integers. In this case $K_n(x;p,a)$
has
$a+1-(-1)-1=a+1>n$ zeros. In both cases,
the number of zeros is greater than the degree of the polynomial and
we have a contradiction. This means $K_n(r+1;p,a)\neq0.$

\medskip \noindent
The proof that $K_n(r-1;p,a)\neq0$ is analogous. Now (\ref{51})
implies that $\displaystyle{K_n(r+1;p,a)=- \frac {C(r)}{A(r)}K_n(r-1;p,a)}$  and clearly
$K_n(r+1;p,a)$ and $K_n(r-1;p,a)$ differ in sign.
\endproof

\begin{Theorem}\label{Krawtchouk} If $a$ is any real number, $a>n-1$,
$n\in\nn$ and  $0<p<1$, the zeros of $K_n(x;p,a)$ lie in the open interval $(0,a).$\end{Theorem}

\proof. Let $a>n-1$, $0<p<1$ and let $n$ and $N$ be integers, such
that $N=\lceil a \rceil$ where $\lceil a \rceil$ denotes the least integer larger than or equal to $a$ .  In the sequence
\begin{equation}K_n(0;p,a),K_n(1;p,a),......,K_n(N;p,a)\label{7}\end{equation}
each term can be considered as a polynomial function of the
parameter $a$ with $ 0<p<1 $ fixed. When a numerical value is
assigned to $a$, we denote the number of variations in sign in the
resulting sequence by $V(a).$ We want to determine $V(a)$ for
$N-1<a\leq N.$\medskip \noindent

When $a=N$, it follows from (\ref{intmass}) that the sequence of
polynomials in (\ref{7}) will have $n$ sign changes since
$K_n(x;p,N)$ is orthogonal for $0<p<1$. This means that $V(N)=n.$

\medskip \noindent If $a$ is assigned any value in the interval $(N-1,N],$ then Lemma
\ref{l} implies that in the resulting sequence
\[K_n(0;p,a),K_n(1;p,a),......,K_n(N;p,a)\]
no two consecutive terms are zero and also that if $K_n(i;p,a)=0$
for $i=1,2,\dots N-1$, then the two adjacent terms:
$K_n(i-1;p,a)$ and $K_n(i+1;p,a)$ differ in sign. Moreover, it
follows directly from (\ref{Kraw2}) that the first term
\begin{equation}\label{K}K_n(0;p,a)=(-a)_n\end{equation} can
never be zero for $a$ in the interval $(N-1,N].$ The last term does
not change sign on $(N-1,N]$ since, by (\ref{77}), \bean K_n(N;p,a)&=&(-a)_n(1-\frac
1p)^n\sum_{i=0}^{n}\frac{(-n)_i(-a+N)_i}{(-a)_i(1-p)^i{i!}}>0~\mbox{for
all}~a\in (N-1,N].\eean
These conditions are sufficient to ensure that the sequence
(\ref{7}) forms a generalised Sturmian sequence (cf. \cite{Porter})
and therefore $V(a)$ remains constant as $a$ increases through the
interval $(N-1,N]$. Hence $V(a)=n~ \mbox{for all}~\\ a \in
(N-1,N].$

\medskip \noindent
Thus for $a>n-1$, $K_n(x;p,a)$ changes
sign $n$ times for $x$ in $(0,N)$ and since the degree
is $n$ we conclude that $K_n(x;p,a)$ has $n$ distinct roots in
$(0,N).$

\medskip \noindent If $r$ is a root of $K_n(x;p,a)$, then $0<r<N$ and it follows from the symmetry relation (\ref{symmetry}) that $a-r$ will be a zero of $K_n(x;1-p,a)$  with
$0<a-r<N$, i.e. $r<a$. We conclude that the zeros of $K_n(x;p,a)$ are in the open
interval $(0,a).$

\endproof

\section{Quasi-orthogonality of $K_n(x;p,a)$, $p<0$ or $p>1$ and $a>0$}
\label{44}

We say that a polynomial $P_n$ of exact degree $n\ge r,$ is discrete
quasi-orthogonal of order $r$ on $[a,b]$ with respect to a weight
function $w(x)$, if (cf.  \cite[p 59]{Brezinsky})
$$\sum_{x=0}^{\infty}x^jP_n(x)w(x)\begin{cases}=0,~\mbox{for} ~j=0,1,\dots,n-r-1\\\neq0,~\mbox{for} ~j=n-r.\end{cases}$$
A more general definition of quasi-orthogonality is given in
(\cite[p 64]{Chihara78}). The Meixner polynomials
$K_n(x;p,a)$ are orthogonal on $(0,\infty)$ for $p<0$, $a<0$ and
as $a$ increases above $0$, the zeros of $K_n(x;p,a)$ depart from
the interval of orthogonality $(0,\infty)$. We prove the
quasi-orthogonality of these polynomials in the following theorem.
\begin{Theorem}\label{Quasi}The polynomials $K_n(x;p,a+k)$ with
$p<0$, $-1<a<0$ and $k=1,2,\dots,n-1$ are quasi-orthogonal of order
$k$ with respect to the weight function $\displaystyle
\frac{p^x(-a)_x}{(p-1)^xx!}$ on $(0,\infty)$.\end{Theorem}
 \proof. The recurrence relation (cf. \cite{Rain}, p 71, Eq. 2)
\begin{eqnarray}K_n(x;p,a)=nK_{n-1}(x;p,a)+K_n(x;p,a+1)\label{112}\end{eqnarray}
shows that $K_n(x;p,a+k)$ can be expressed as a linear combination
of $K_n(x;p,a),K_{n-1}(x;p,a),\dots,K_{n-k}(x;p,a)$ and, since $a<0,$
it follows from (\ref{1}) and (\ref{m}) that
\[\sum_{x=0}^{\infty}x^jK_n(x;p,a+k)\frac{p^x(-a)_x}{(p-1)^xx!}= 0~\mbox{for}~j=0,1,\dots,n-k-1.\]\endproof
\begin{Remark}By a change of variable, the result in Theorem
\ref{Quasi} can be written as that the polynomials $K_n(x;p,a)$
are quasi-orthogonal of order $k$ on $(0,\infty)$, for $p<0$ and $k-1<a<k$,
$k=1,2,\dots,n-1$ with respect to the weight function $\displaystyle
\frac{p^x(-a+k)_x}{(p-1)^xx!}$.\end{Remark}
The zeros of quasi-orthogonal polynomials are not necessarily all in
the interval of orthogonality, but we can say the following from
\cite[Thm 2]{Brezinsky}.
\begin{Corollary} The Meixner polynomials $K_n(x;p,a)$, with $p<0,~ k-1<a<k$ have at least
$n-k$ zeros in $(0,\infty)$ when $k=1,2,\dots,n-1.$
\end{Corollary}


\medskip \noindent
In order to specify the location of the remaining single zero of $ K_n(x;p,a+1),p<0,-1<a<0$ where we have quasi-orthogonality of order $1$,
we consider the monic Krawtchouk polynomials $\tilde{K}_n(x;p,a)=p^n K_n(x;p,a)$.

\begin{Theorem} If $p<0$ and $-1<a<0$, then the smallest zero of $K_n(x;p,a+1)$ (or equivalently $ \tilde{K}_n(x;p,a+1)$) is negative.
\end{Theorem}
\proof. The recurrence relation (\ref{112}) can be written as
\begin{equation}\tilde{K}_n(x;p,a+1)= \tilde{K}_n(x;p,a)-np \tilde{K}_{n-1}(x;p,a)\label{TTRR*}\end{equation} and according to \cite[Thm 4]{Joulak} we have to show that $\displaystyle~np<\frac{\tilde{K}_{n}(0,p,a)}{\tilde{K}_{n-1}(0,p,a)}<0$ which follows immediately from (\ref{K}). \endproof

Joulak's results (cf. \cite[Thms 8, 9]{Joulak}) also gives some information about the location of the zeros when we have quasi-orthogonality of order $2$.

\begin{Theorem} If $p<0$, $-1<a<0$ and $n>(a+2)(1-p)$ then all the zeros of $K_n(x;p,a+2)$ are nonnegative and simple.
\end{Theorem}
\proof.
Iterating (\ref{TTRR*}) we obtain
\[\tilde{K}_n(x;p,a+2)=\tilde{K}_n(x;p,a)-2np\tilde{K}_{n-1}(x;p,a)+b_n\tilde{K}_{n-2}(x;p,a)\] where $b_n=n(n-1)p^2$. Letting $a=N$ and replacing $n$ by $n-1$ in (\ref{ttR}) yields
\[\tilde{K}_n(x;p;a)=(x-B_{n-1})\tilde{K}_{n-1}(x;p;a)-pC_{n-1}\tilde{K}_{n-2}(x;p,a).\label{100}\]
From \cite[Thm 8]{Joulak} all the zeros of $\tilde{K}_n(x;p,a+2)$ are real and simple if $b_n<pC_{n-1}$ which gives the condition
\[ n>(a+2)(1-p).\]
Furthermore, the smallest zero (and hence all of the zeros) of $K_{n}(x;p,a+2)$ is nonnegative if and only if (cf. \cite[Thm 9]{Joulak})
\[ \frac{\tilde{K}_{n}(0,p,a)}{\tilde{K}_{n-2}(0,p,a)}-2np\frac{\tilde{K}_{n-1}(0,p,a)}{\tilde{K}_{n-2}(0,p,a)}+n(n-1)p^2\geq 0.\]
It follows from (\ref{K}) that the left-hand side simplifies to $p^2(a+1)(a+2)$, which is positive by the assumptions. This completes the proof.
\endproof

Note that the sufficient condition for the positivity of the zeros from \cite[Thm 4]{Brezinsky} does not apply here, since $b_n>0$.

\medskip \noindent
Analogous results can be obtained for the polynomials
$K_n(x;p,a)$,  $p>1,~ k-1<a<k$.

\section{The zeros of $K_{n}(x;p,N)$, $0<p<1$ and $n=N+1,N+2,\dots$}\label{3}
We now turn our attention to the zeros of the polynomials $K_{n}(x;p,N)$, $0<p<1$
and $n,~N\in \nn$  for degrees higher than what follows from
orthogonality, i.e. for $n= N+k, k\in\nn$.
The case when $n=N+1$ was first studied by J Sylvester (cf. \cite{Sylvester}) using
Sylvester type determinants and more recently by R Askey and G Wilson (cf. \cite{Ask}).
The connection of orthogonal polynomials
with tridiagonal matrices whose entries come from the recurrence
coefficients of discrete orthogonal polynomials is made explicit in
\cite{AskSylv} and \cite{Holtz}. Both these papers show that
$K_n(x;p,N)$ has zeros at $0,1,\dots,N$ when $n=N+1$, $n,~N\in \nn$
and $0<p<1$. For completeness and the convenience
of the reader, we provide a direct proof.

\begin{Lemma} (cf. \cite{Holtz})\label{Holtz} Let $0<p<1$ and $N$ a positive integer,
the polynomials $K_{N+1}(x;p,N)=(x)(x-1)\dots(x-N)(\frac{1}{p})^{N+1}$ and have $N+1$ real zeros $x=0,1,\dots,N.$\end{Lemma}

\proof.
Since
\[(-N+k)_{N+1-k}\begin{cases}=0~\mbox{for}~k=1,\dots,N\\=1~\mbox{for}~k=N+1,\end{cases}\]
it follows from (\ref{Kraw2}) that \begin{eqnarray*}K_{N+1}(x;p,N)&=&\frac{(-N-1)_{N+1}(-x)_{N+1}}{p^{N+1}(N+1)!}\\
&=&\frac{(x)(x-1)\dots(x-N)}{p^{N+1}}.\end{eqnarray*}
\endproof

\begin{Corollary}\label{n=N+2} For $0<p<1$ and $N$ a positive integer,
the polynomial $K_{N+2}(x;p,N)$
has $N+2$ real zeros \bean x=0,1,\dots,N, N+1-p(N+2).\eean

\end{Corollary}\vspace{0.2cm}

\proof.
Letting $n=N+2$ in (\ref{product}) we obtain $\dis{K_{N+2}(x;p,N)=
K_{N+1}(x,p,N)\left(N+2+\frac{x-(N+1)}{p}\right)}$ which yields the stated result.
\endproof

\begin{Corollary}\label{n=N+3} For $0<p<1$ and $N$ a positive integer, the polynomial $K_{N+3}(x;p,N)$
has at least $N+1$ real zeros $x=0,1,\dots,N$. Furthermore, the
remaining two zeros will be real and distinct when
\begin{eqnarray*}0< p<\frac{1}{2}\left(1-\sqrt{\frac{N+2}{N+3}}\right)&~\mbox{or}~&
\frac{1}{2}\left(1+\sqrt{\frac{N+2}{N+3}}\right)<p<1.\end{eqnarray*}
\end{Corollary}\vspace{0.2cm}

\proof. It follows from (\ref{product}) that $ K_{N+3}(x;p,N) = K_{N+1}(x;p,N)~p_2(x)$, where \bean
p_2(x)
&=&\frac{1}{p^2}\left(x^2+ (6p+2 Np-3-2 N) x +(N+3)(N+2)p^2 -
2 (N+1) (N+3)p + (N+1) (N+2)\right).\eean The zeros of the quadratic $p_2(x)$ are real and distinct when the discriminant
$\displaystyle{\frac{1}{p^4}\left(1-4p(N+3)+4p^2(N+3)\right)}$ is
positive, i.e. for \bean p<\frac{1}{2}\left(1-\sqrt{\frac{N+2}{N+3}}\right) \mbox{  or  }
p>\frac{1}{2}\left(1+\sqrt{\frac{N+2}{N+3}}\right). \eean\endproof

\noindent It is difficult to determine the exact location of the zeros in
the general case $K_n(x;p,N)$, $n=N+k$, $k\in \nn$. There are only $N+1$ or $N+2$ real zeros and
the lack of orthogonality makes potential theory difficult to apply.
Naturally, one can use (\ref{product}) to consider the polynomials
$K_n(x;p,-N)$ for $0<p<1$, $N\in \nn$ instead. In general, it
suffices to investigate the zeros of $K_n(x,p,a)$ for $a<0$ when
$0<p\leq 1/2$ (or $1/2\leq p<1$), because it follows from the
symmetry relation (\ref{symmetry}) that if $x$ is a zero of
$K_n(x;p,a)$ then $a-x$ is a zero of $K_n(x;1-p,a)$. Taking also
into consideration the complex conjugate pairs, geometrically it
means that the zeros of $K_n(x;1-p,a)$ are the mirror image of the
zeros of $K_n(x;p,a)$ with respect to the axis $\Re x=a/2$ when
$a<0$ and $0<p<1$.

\begin{figure}[h]
\centering
\parbox{6cm}{\includegraphics[height=4cm]{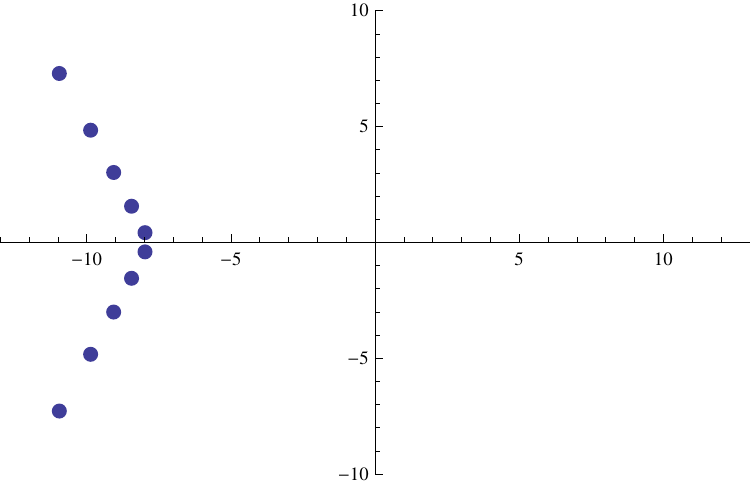}
\caption{The zeros of $K_{10}(x;0.8,-8.2)$.}
\label{fig2}}
\qquad
\parbox{6cm}{
\includegraphics[height=4cm]{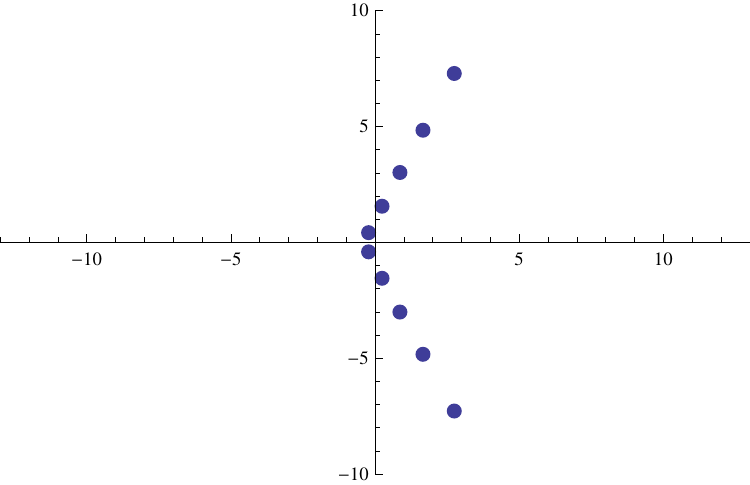}
\caption{The zeros of $K_{10}(x;0.2,-8.2)$.}
\label{fig3}}
\centering
\parbox{6cm}{
\includegraphics[height=4cm]{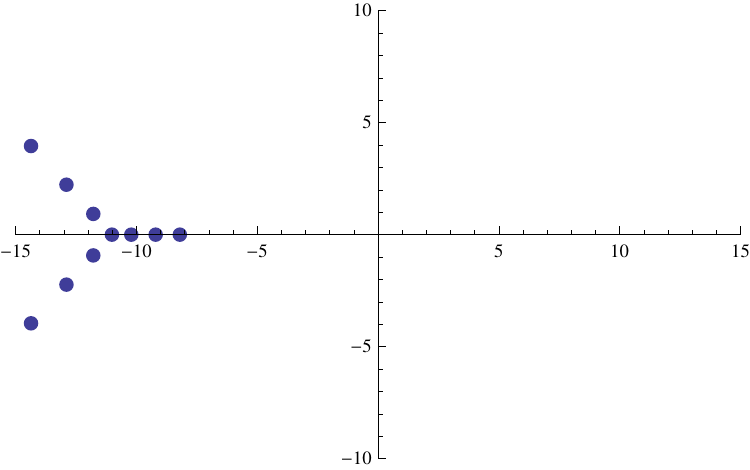}
\caption{The zeros of $K_{10}(x;0.94,-8.2)$.}
\label{fig4}}
\qquad
\parbox{6cm}{\includegraphics[height=4cm]{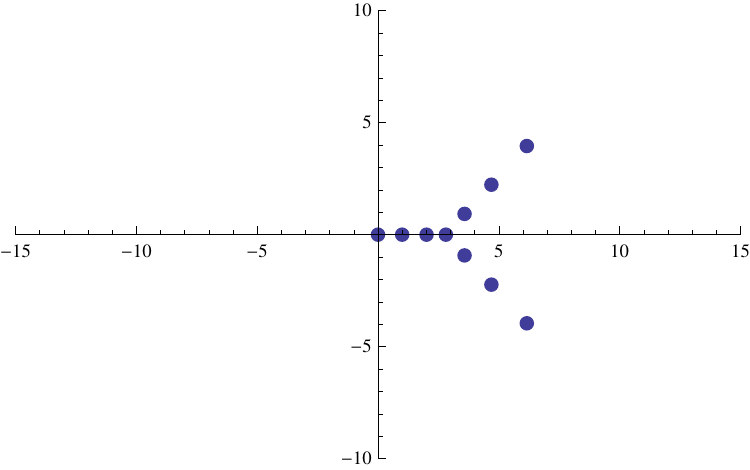}
\caption{The zeros of $K_{10}(x;0.06,-8.2)$.}
\label{fig5}}
\end{figure}

\medskip\noindent Figures \ref{fig2}, \ref{fig3}, \ref{fig4} and \ref{fig5} show the zeros of $K_n(x;p,a)$ when $n=10$, $a=-8.2$ for different values of $p$, clearly illustrating the symmetry with respect to $\Re x=-4.1$.

\begin{figure}[h!]
\begin{center}
\includegraphics[height=4cm]{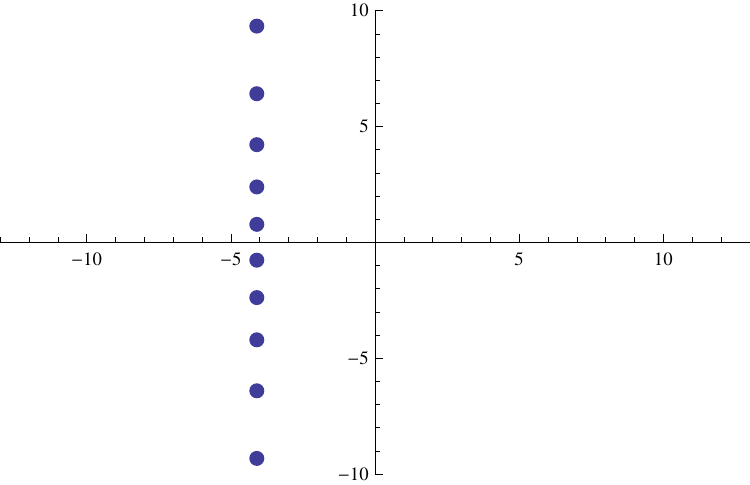}\caption{The zeros of $K_{10}(x;0.5,-8.2)$}
\label{fi1}
\end{center}
\end{figure}

\medskip\noindent
The numerical examples show that the zeros of polynomials $K_n(x;p,a)$, $a<0$, $0<p<1$  seem to lie on rays starting
from the $x$ axis. At the special case $p=1/2$ all the zeros lie on the line $\Re x=a/2$. This special case
when $p=\frac 12$ is illustrated in Figure \ref{fi1}.

\begin{Remark}
Letting $a=-n$, $b=-x$, $c=-N-1$ and $z=\frac 1p$ in the contiguous relation (cf. \cite[Section 33]{Rain})
\[ [a+(b-c)z] \phantom{ }~_2F_1(a,b;c;z)=a(1-z)\phantom{ }~_2F_1(a+1,b;c;z)-\frac{(c-a)(c-b)z}{c}\phantom{ }~_2F_1(a,b;c+1;z) \]
we obtain
\[ (x-N-1)K_n(x;p,N)=n(n-N-2)(1-p)K_{n-1}(x;p,N+1)+(x+np-N-1)K_n(x;p,N+1). \]
Using this it is easy to show by induction that
\[ (x-N-1)(x-N-2)\dots (x-N-j)K_n(x;p,N)=\sum_{i=0}^j q_i(x)K_{n-i} (x;p,N+j), \]
where $q_i(x)$ are polynomials of degree $j-i$. From this it follows (cf.  \cite[p. 160]{Brezinsky}) that $(x-N-1)(x-N-2)\dots (x-N-k)K_n(x;p,N)$ (with $n=N+k$) is quasi-orthogonal of degree $N+2k$ and order $2k$. However, unlike in Section \ref{44} this does not lead to new information about real zeros, since these are already guaranteed by formula (\ref{product}) and Lemma \ref{Holtz}.
\end{Remark}

\section{The zeros of $K_{n}(x;p,a)$, $p\to 0$ and $a\in \rr$}

Next, we consider $K_n(x;p,a)$ when $a\in \rr$  and prove that when
$p\rightarrow 0$ and $a\neq0$ all the zeros of $K_n(x;p,a)$,
$n=1,2,\dots$ approach non-negative integer values.
Note that this theorem holds for any $a\in\rr$ which implies that when $a=N$,
$N\in \nn$ and $p>0$, the zeros of the Krawtchouk polynomials approach the mass points of the weight function as $p\to 0.$

\begin{Theorem}\label{tinfinity} The $n$ zeros of the polynomial \[K_n(x;p,a),~a\in\rr
,~n=0,1,\dots\] approach the points $x=0,~1,\dots,~n-1$ when
$p\rightarrow 0$.\end{Theorem}

\proof.
\bean K_n(x;p,a)&=&(a)_n+
\dots+\frac{(a+n-1)(-n)_{n-1}(-x)_{n-1}}{(n-1)!p^{n-1}}+\frac{(x)(x-1)\dots(x-n+1)}{p^n}.\eean
For any $n=0,1,2,\dots$ the function\begin{eqnarray*}
p^nK_n(x;p,a)=p^n(a)_n +
\dots+\frac{p(a+n-1)(-n)_{n-1}(-x)_{n-1}}{(n-1)!}+(x)(x-1)\dots(x-n+1)\end{eqnarray*}
regarded as an $n$th degree polynomial in $x$ with real parameters
$a$ and $p$ has the same zeros as $K_n(x;p,a)$. Since
\bean\lim_{p\to 0}p^nK_n(x;p,a)&=&(x)(x-1)\dots(x-n+1),\eean the zeros of
$\displaystyle{p^nK_n(x;p,a)}$ and hence the zeros of $K_n(x;p,a)$ tend to the zeros of $x(x-1)(x-2)\dots(x-n+1),$ which is to say $x=0,1,2,\dots,n-1.$\endproof

\noindent This Theorem implies that for sufficiently small $p$ all the zeros of $K_n(x;p,a)$ are real. An analogous result can be proved for Charlier polynomials.

\medskip
\noindent{\bf Acknowledgement}
The first two authors would like to thank the Institute for Biomathematics and Biometry at the Helmholtz Zentrum M\"{u}nchen for providing accommodation during a research visit where this a research was completed.

\end{document}